\documentclass[11pt]{amsart}
\usepackage{amsmath,amsthm, amscd, amssymb, amsfonts}

\newcommand{\ku}{\Bbbk}

\newcommand\toba{{\mathfrak B }}

\newcommand\tobar{ {\widehat{\mathfrak B}}_r  }
\newcommand\tobados{ {\widehat{\mathfrak B}}_2  }

\newcommand{\trid}{\triangleright}

\newcommand{\g}{{\mathfrak g}}

\newcommand{\Z}{{\mathbb Z}}

\newcommand{\Cc}{{\mathbb C}}

\newcommand{\otravez}{{1\le i \le \theta}}
\newcommand{\otraavez}{{1\le i \neq j \le \theta}}
\newcommand{\otrvez}{{1\le i,j \le \theta}}

\numberwithin{equation}{subsection}

\newcommand{\Aut}{\operatorname{Aut}}

\newcommand{\End}{\operatorname{End}}

\newcommand{\ord}{\operatorname{ord}}

\newcommand{\ev}{\operatorname{ev}}

\newcommand{\Ext}{\operatorname{Ext}}

\newcommand{\gd}{\operatorname{gldim}}

\newcommand{\id}{\mathop{\rm id\,}}

\theoremstyle{plain}

\newtheorem{theorem}{Theorem}[subsection]
\newtheorem{lema}[theorem]{Lemma}

\newtheorem{prop}[theorem]{Proposition}

\theoremstyle{definition}
\newtheorem{definition}[theorem]{Definition}

\theoremstyle{remark}
\newtheorem{obs}[theorem]{Remark}

\def\pf{\begin{proof}}
\def\epf{\end{proof}}

\theoremstyle{remark}

\begin{document}

\renewcommand{\baselinestretch}{1.2}
\thispagestyle{empty}
\title{Some remarks on Nichols algebras}
\author[Andruskiewitsch]{ Nicol\'as Andruskiewitsch}
\address{Facultad de Matem\'atica, Astronom\'\i a y F\'\i sica,
Universidad Nacional de C\'ordoba
\newline
\indent CIEM -- CONICET
\newline
\indent (5000) Ciudad Universitaria, C\'ordoba, Argentina}
\email{andrus@mate.uncor.edu}
\thanks{This work was partially supported by CONICET,
Agencia C\'ordoba Ciencia, ANPCyT  and Secyt (UNC)}
\date{\today}
\begin{abstract} Two algebras can be attached to
a braided vector space $(V, c)$
in an intrinsic way; the FRT-bialgebra and the Nichols algebra
$\toba(V, c)$. The FRT-bialgebra plays the r\^ole of the algebra
of quantum matrices, whereas the r\^ole of the Nichols algebra is
less understood. Some authors call $\toba(V)$ a quantum symmetric
algebra. The purpose of this paper is to discuss some properties
of certain Nichols algebras, in an attempt to establish  classes
of Nichols algebras which are worth of further study.
\end{abstract}
\maketitle

\subsection{Definitions and examples of Nichols algebras}
Let $(V, c)$ be a braided vector space, that is, $V$ is a
finite-dimensional complex vector space and  $c: V \otimes V \to V
\otimes V$ is an invertible solution of the braid equation:
$(c\otimes \text{id})(\text{id} \otimes c)(c\otimes \text{id}) =
(\text{id} \otimes c)(c\otimes \text{id}) (\text{id} \otimes c)$.
There is a remarkable braided graded Hopf algebra $\mathfrak B (V,
c) = \oplus_{n \ge 0} \mathfrak B^n(V, c)$, which is connected,
generated in degree one, with $\mathfrak B^1(V) \simeq V$ as
braided vector spaces,  such that all its primitive elements have
degree one; and which is unique with respect to these properties.
Algebras of this kind
appeared naturally in our approach to classification of
pointed Hopf algebras \cite{AS1, AS-cambr} but we quickly realized
they were already known to several authors under various
presentations. We first briefly recall different definitions of
Nichols algebras and survey examples of classes of Nichols
algebras that are known. A detailed exposition on Nichols algebras
can be found in \cite{AS-cambr}.

\bigbreak We shall simply write $\mathfrak B (V) = \mathfrak
B (V, c)$ omitting the reference to $c$ unless it is needed.
We shall always assume that
$c$ is rigid, \emph{i.~e.} the associated map $c^{\flat}: V^* \otimes V\to V \otimes V^*$
is also bijective (this is the case in all the examples below).
Here $c^{\flat} = (\ev_V \otimes \id_{V\otimes V^*})(\id_{V^*}\otimes c\otimes \id_{V^*})
(\id_{V^*\otimes V}\otimes \ev_V^*)$.

\bigbreak This remarkable braided Hopf algebra was first
described by W. Nichols in his thesis \cite{N}, as the invariant
part of his "bialgebras of type one". In his honor, $\mathfrak B
(V)$ is called the Nichols algebra of the braided vector space
$(V, c)$. There are several ways to present $\mathfrak B (V)$.

\bigbreak
Consider $T(V) \otimes T(V)$ as an algebra with the product 'twisted'
by $c$. Then $T(V)$ is a braided Hopf algebra, with the comultiplication
uniquely defined by $\Delta(v) = v\otimes 1 + 1 \otimes v$, $v\in V$.
Let $I(V)$ be the largest Hopf ideal generated by homogeneous elements
of degree greater than 1; then $\toba (V) := T(V)/I(V)$ satisfies all
the properties listed above \cite[Prop. 2.2]{AS-cambr}.

\bigbreak The vector space $T(V)$ has another structure of
coalgebra, the free coalgebra over $V$; let us denote it by
$t(V)$. M. Rosso observed that it admits a 'quantum shuffle
product', so that $t(V)$ is also a braided Hopf algebra, called
the quantum shuffle algebra. The canonical map $\Omega: T(V) \to
t(V)$ turns out to be a map of braided Hopf algebras; the image of
$\Omega$, that is the subalgebra of $t(V)$ generated by $V$, is
the Nichols algebra of $V$. The nilpotent part $U_q^+(\mathfrak
g)$ of a quantized enveloping algebra was characterized in this
way by Rosso as the Nichols algebra of a suitable braided vector space
$\mathfrak h$ \cite{Ro1, Ro2}. Results in the same spirit were
also obtained by  J. A. Green \cite{Gr}.

\bigbreak
Now, the components of the graded map $\Omega$, that is $\Omega^{(m)}:
T^m(V) \to T^m(V)$ are the so-called ``quantum symmetrizers''
defined through the action of the braid group on $T^m(V)$.
Therefore, the Nichols algebra of $(V, c)$ coincides with the quantum exterior
algebra of $(V, -c)$, defined by   S. L. Woronowicz \cite{W}. Indeed, the quantum symmetrizers
of $-c$ are the quantum antisymmetrizers of $c$.

\bigbreak G. Lusztig characterized $U_q^+(\mathfrak g)$ as the
quotient of a $T(\mathfrak h)$ by the radical of an invariant
bilinear form \cite{L3}. This is indeed a general fact; the ideal
$I(V)$ is always the radical of an invariant bilinear form
\cite{AG}.

\bigbreak M. Rosso found that the Nichols algebra $\toba(V)$
of a braided vector space of diagonal type has
always a "PBW-basis" in terms of the so-called Lyndon words
\cite{Ro3}. Related work was done by  V. K. Kharchenko, who also
studied abstractly Nichols algebras from various points of view
\cite{Kh1, Kh2, Kh3}.

\bigbreak The following two questions arise from classification
problems of Hopf algebras \cite{AS1, AS-cambr}. Answers to both
questions are needed to classify Hopf algebras of certain types.

\bigbreak
\begin{itemize}
\item[$\circ$] Under which conditions
on $(V, c)$ is $\mathfrak B(V)$ finite-dimensional, respectively
of finite Gelfand-Kirillov dimension?

\bigbreak
\item[$\circ$] For those pairs with a positive answer to the preceding question,
give an explicit presentation of $\mathfrak B(V)$; that is, find a minimal
set of generators of the ideal $I(V)$.

\end{itemize}

\bigbreak
The study of $\mathfrak B(V)$ is very difficult;
neither the subalgebra of the quantum shuffle algebra generated by $V$,
nor the Lyndon words, nor the ideal $I(V)$ have an explicit description.

\bigbreak
\begin{itemize}
\item[$\circ$] In particular, we do not know if the ideal
$I(V)$ is finitely generated. \end{itemize}

\bigbreak There is little hope to perform explicit computations
with a computer program without a positive answer to this
question. However, let $\tobar (V) = T(V) / J_r$, where $J_r$ is
the two-sided ideal generated by the kernels of $\Omega^{(m)}$,
$n\le r$; these are braided Hopf algebras and we have epimorphisms
$\tobar (V) \to \toba (V)$ for all $r\ge 2$. Hence, if one of the
algebras $\tobar (V)$ is finite-dimensional, or has finite
Gelfand-Kirillov dimension, so does $\toba (V)$. In the first
case, under favorable hypothesis we may conclude that $\tobar (V)
\simeq \toba (V)$, see \cite[Th. 6.4]{AG2}.

\bigbreak
There are several classes of braided vector spaces which seem
to be of special interest.

\begin{itemize}

\bigbreak
\item We say that $(V, c)$ is \emph{of diagonal type} if there exists
a basis $x_1, \dots, x_{\theta}$ of $V$, and non-zero scalars
$q_{ij}$ such that $c(x_i\otimes x_j) = q_{ij} x_j \otimes x_i$,
$\otrvez$. \end{itemize}

 \bigbreak
Nichols algebras of these braided vector spaces appear
naturally in the classification of pointed Hopf algebras with
abelian coradical, and also in the theory of quantum groups.

Namely, let $(a_{ij})_{\otrvez}$ be a generalized Cartan matrix;
let $\mathfrak h$ be a vector space with a a basis $x_1, \dots,
x_{\theta}$, let $q$ be a non-zero scalar and let $c$ be given by
\begin{equation}\label{quantum}
c(x_i\otimes x_j) = q^{a_{ij}} x_j \otimes x_i, \qquad \otrvez.
\end{equation}
 Then $\toba(\mathfrak h) =
U_q^+(\mathfrak g)$ if $q$ is not a root of 1 \cite{L3, Ro1, Ro2}, and
$\toba(\mathfrak h) = \mathfrak u_q^+(\mathfrak g)$ if $q\neq 1$ is a root of 1 (under some hypothesis
on the order of $q$) \cite{Ro1, Ro2, Mu}.

 \bigbreak
 \begin{itemize}
\item We say that $(V, c)$ is \emph{of rack type} if there exists a basis $X$ of $V$,
a function $\trid: X\times X \to X$ and non-zero scalars $q_{ij}$
such that
\begin{equation}\label{rack}
c(i\otimes j) = q_{ij} \, i \trid j \otimes i,
\end{equation}
$i,j \in X$. Then $(X, \trid)$ is a rack and $q_{ij}$ is a rack 2-cocycle
with coefficients in $\ku^{\times}$, see for example \cite{G-cont, AG2}.
\end{itemize}

These
braided vector spaces appear naturally in the classification of
pointed Hopf algebras.

\bigbreak
\begin{itemize}
\item We say that $(V, c)$ is \emph{of Jordanian type} if there exists
a basis $x_1, \dots, x_{\theta}$ of $V$, and a non-zero scalar $q$
such that $c(x_i\otimes x_1) = q x_1 \otimes x_i$, $c(x_i\otimes
x_j) = (q \, x_j + x_{j-1}) \otimes x_i$, $\otravez$, $2\le j \le
\theta$. \end{itemize}

These braided vector spaces appear  in the
classification of pointed Hopf algebras with  coradical $\Z$.

\bigbreak
\begin{itemize}
\item We say that $(V, c)$ is \emph{of Hecke type}
if $(c - q)(c+1) = 0$, for some non-zero scalar $q$.

\bigbreak
\item We say that $(V, c)$ is \emph{of quantum group type} if $V$
is a module over some quantized enveloping algebra $U_q(\g)$ and $c$ arises
from the action of the corresponding universal $R$-matrix.
\end{itemize}

\bigbreak
Here is what is known about the problems stated above.

\bigbreak
\begin{itemize}
\item $\mathfrak B(V) = T(V)$ generically.
That is, consider the locally closed space of all $c\in
\End(V\otimes V)$ which are invertible solutions of the braid
equation; then the subset of those $c$ such that $\mathfrak B(V,
c) = T(V)$ contains a non-empty open subset.

\bigbreak
\item Assume that $(V, c)$ is of diagonal type, where the $q_{ii}$'s are
positive and different from one, $\otravez$. Then $\toba(V)$ has
finite Gelfand-Kirillov dimension if and only if $q_{ij}q_{ji} =
q_{ii}^{a_{ij}}$ for some Cartan matrix of finite type \cite{Ro2}.

\bigbreak
\item Assume that $(V, c)$ is of diagonal type, that the
$q_{ii}$'s are are roots of 1 but not 1,
and that $q_{ij}q_{ji} = q_{ii}^{a_{ij}}$, $\otrvez$, where
$a_{ij} \in \Z$, $a_{ii} = 2$, $\ord q_{ii} < a_{ij} \le 0$ if $\otraavez$.
Then $(a_{ij})_{\otrvez}$ is a generalized Cartan matrix and, under
suitable conditions, $\toba(V)$ has finite dimension if and only
if $(a_{ij})_{\otrvez}$ is of finite type \cite{AS-adv}.
\end{itemize}

\bigbreak In these two cases, the calculation of $\toba(V)$ is
reduced to the calculation of $\toba(\mathfrak h)$ as in
\eqref{quantum}; but the last requires deep facts on
representation theory, and the action of the quantum Weyl group
defined by Lusztig.

\bigbreak\begin{itemize}
\item There are a few examples of $(V, c)$ of diagonal type with
finite-dimensional $\toba(V)$, due to Nichols \cite{N} and Gra\~na
\cite{G-cont}, besides those in the last item (they are listed in
\cite[Section 3.3]{AS-cambr}).

\bigbreak
\item If $c$ is of Hecke type and $q$ is not a root of 1, or if
$q=1$, then the Nichols algebra is quadratic: $\toba(V) = \tobados
(V) $. Furthermore, the quadratic dual is also a Nichols algebra:
$\toba(V)^{!} = \toba(V^*, -q^{-1}c^t)$; see Proposition \ref{cuadratic} below.

\bigbreak
\item If $c$ is of diagonal type,
information about $\det \Omega^{(m)}$ is
given in \cite{FG}.

\bigbreak
\item There are a few examples of $(V, c)$ of rack type with
finite-dimensional $\toba(V)$ \cite{MS, G-cont, AG2, G-zoo}. See
the table in the Appendix for a flavor of the kind of algebras
obtained.

\bigbreak
\item[$\circ$] Almost nothing is known about Nichols algebras of quantum group type,
except when they are of Hecke type. To begin with, it would be interesting
to know what happens when $V = L(n)$ is a highest weight module over $U_q(\mathfrak{sl} (2))$,
$n \ge 3$ (if $n = 1$ it is of Hecke type, if $n = 2$ it seems to be known but I do not have a
reference).

\bigbreak
\item[$\circ$] Nichols algebras of Jordanian type were not considered in the literature,
to my knowledge. It is likely that the quantum Jordanian plane is a Nichols algebra
of Jordanian type.
\end{itemize}

\bigbreak In conclusion, a Nichols algebra may be
finite-dimensional or not, have finite Gelfand-Kirillov dimension
or not, and there is no general technique, up to now, to
explicitly decide for a given braided vector space, what is the
case for its Nichols algebra\footnote{Some techniques are
available for specific classes of Nichols algebras, \emph{e.~g.}
for Nichols algebras of diagonal type, as already said.}.

This indicates that Nichols algebras do not have to be studied
through a general approach, but splitting the category of
braided vector spaces in classes.

\bigbreak
\begin{obs} Another important question in the classification
of Hopf algebras is the following. Let $B  = \oplus_{n \ge 0} B^n$
be a braided graded Hopf algebra, connected, generated in degree
one, and denote by $V$ the braided vector subspace $B^1$ of $B$.
Is it possible to conclude that $B$ is the Nichols algebra of $V$,
\emph{i.~e.} that all its primitive elements have degree one,
under some abstract conditions? Partial positive answers to this
question are given in \cite[Th. 7.6]{AS4} (finite-dimensional
case) and \cite[Lemma 5.1]{AS6} (finite Gelfand-Kirillov dimension
case).

\end{obs}

\subsection{Some properties of some Nichols algebras}
It is  clear that no fine ring theoretical properties can be
established for Nichols algebras in general. But
there might be a suitable class of
braided vector spaces whose Nichols algebras  desserve attention
from this point of view.

\subsubsection{Graded algebras}

We shall only consider graded algebras $R = \oplus_{i\ge 0} R_i$,
which are finitely generated and $R_0 = \Cc$.

\begin{definition}
A graded algebra $R$ is {\it AS-regular} if it has finite global
dimension $d$, finite Gelfand-Kirillov dimension and  is
AS-Gorenstein. Thus, $ \Ext^i_R(\Cc, R) = 0$, if $i\neq d$,  and
$= \Cc$ if $i = d$.
\end{definition}

This class of graded algebras has been intensively investigated in the
last years; AS is in honor of Artin and Schelter.
The study of the category of graded modules of such an algebra
has a strong geometrical flavor; this is usually
called a noncommutative projective space.
In particular, the space of all "point modules" is a genuine projective space
which provides important information on the full category.
See \cite{S, SV}.
The homological conditions are designed to insure good regularity properties;
we refer again to \cite{S, SV} and references therein.

\bigbreak
\begin{itemize}
\item[$\circ$] When is a Nichols algebra AS-regular?

\bigbreak
\item[$\circ$] Let $\toba(V)$ be a Nichols algebra which is a domain with finite Gelfand-Kirillov dimension.
Is it AS-regular?

\end{itemize}

\bigbreak
Some insight about these questions is explained in the next subsections.

\subsubsection{Koszul algebras} We first recall some well-known facts
about Koszul algebras.

A \emph{quadratic algebra} is a graded algebra $A = \oplus_{n \ge
0} A_n$ generated in degree one with relations in degree 2; that
is $A \simeq T(V) / \langle R \rangle$, where $V = A_1$ and  $R
\subset V \otimes V$ is the kernel of the multiplication. We shall
denote $A = (V, R)$. The quadratic dual of a quadratic algebra $A
= (V, R)$ is $A^{!} = (V^*, R^{\bot})$. By \cite{L},
\begin{equation}
A^{!} \simeq E(A) := \text{ the subalgebra of } \Ext^*_A(\Cc, \Cc)
\text{ generated by }\Ext^1_A(\Cc, \Cc).
\end{equation}

\bigbreak A \emph{graded Koszul} algebra is a quadratic algebra
$A$ such that $A^{!} \simeq \Ext^{\bullet}(\Cc, \Cc)$ as graded algebras.

\begin{lema} \label{lemasmith}  Let $A$ be any graded connected algebra.

(a).  \cite[1.4 and 5.9]{s}. If $\gd A < \infty$ then $\dim A^{!}< \infty$. The converse
holds if $A$ is Koszul.

(b).  \cite[5.10]{s}. If $A$ is Koszul and has finite global dimension, then $A$ is
AS-Gorenstein if and only if $A^{!}$ is Frobenius. \qed
\end{lema}

We can now decide when a Nichols algebra of Hecke type is
AS-regular.

\begin{prop}\label{cuadratic}
Let $(V, c)$ be a braided vector space and assume that $c$
satisfies a Hecke-type condition with label $q$, $q = 1$ or not a
root of 1. Then $\toba(V)$ is AS-regular if and only if it has finite Gelfand-Kirillov
dimension  and the dimension of $\toba(V^*)$ is finite.
\end{prop}

\pf By \cite[Prop. 3.3.1]{AA}, see also \cite[Prop.
3.4]{AS-cambr}, the Nichols algebra $\toba(V)$ is quadratic, and
its quadratic dual is $\toba(V)^{!} = \toba(V^*)$, the Nichols
algebra with respect to $-q^{-1}c^t$. By \cite{Gu, Wa}, $\toba(V)$
is Koszul. By Lemma  \ref{lemasmith} (a), if $\toba(V)$ is
AS-regular then $\dim \toba(V^*) < \infty$.
Conversely, assume that $\dim \toba(V^*) < \infty$. By Lemma  \ref{lemasmith}
part (a), $\gd A < \infty$; and by part (b), $\toba(V)$ is
AS-Gorenstein. Indeed, $\toba(V^*)$ is a braided Hopf algebra;
hence it is Frobenius whenever finite-dimensional.
\epf

\bigbreak
\begin{obs}\label{zhang} I do not know if the hypothesis on the Gelfand-Kirillov
dimension can be removed. There are examples of Koszul algebras
where $A^!$ is finite-dimensional and Frobenius but $A$ has
infinite GK-dimension. I am indebted to James Zhang for pointing
out this to me. We do know the Hilbert series of $\toba(V)$:
$H(\toba(V)) (t) = \dfrac{1}{H(\toba(V^*), -q^{-1}c^t) (-t)}$, by
\cite[Th. 2.11.1]{BGS}, and $H(\toba(V^*))(t)$ is a polynomial.
\end{obs}

\bigbreak
\begin{obs}\label{gurka} D. Gurevich studied intensively Nichols algebras
of Hecke type,  provided ways to construct explicit examples, and classified those
such that $H(\toba(V^*), -q^{-1}c^t)$ is  a polynomial of degree two \cite{Gu}.
\end{obs}

\bigbreak
\begin{obs}\label{roos} A quadratic Nichols algebra is not
necessarily Koszul; see \cite{R}.
\end{obs}

\subsubsection{Nichols algebras related to quantum groups}

\begin{prop}\label{quantum-prop}
Let $(a_{ij})_{\otrvez}$ be a Cartan matrix of finite type,
let $q$ be a non-zero scalar and let $(\mathfrak h, c)$ be the braided  vector
space as in \eqref{quantum}. Then

(a). \cite[Th. 4.7]{GL} $\toba(\mathfrak h) = U_q^+(\mathfrak g)$ is AS-regular when
$q$ is not a root of 1.

(b). \cite[Th. 2.5]{GK} $\toba(\mathfrak h)= \mathfrak u_q^+(\mathfrak g)$ is not
AS-regular if $1\neq q$ is a root of 1. \qed
\end{prop}

It would be interesting to have another proof of (a) in the spirit of \cite{deck, GK}.
Namely, to consider the algebra filtration given by the PBW-basis;
the associated graded algebra is a quantum
linear space \cite{deck}, hence a Nichols algebra of Hecke
type; therefore it has finite global dimension. Then lift this information
by a spectral sequence argument. Note also that
Proposition \ref{quantum-prop} extends to the multiparametric case
without difficulties. In view of Propositions \ref{cuadratic} and \ref{quantum-prop},
it is tempting to suggest the following questions.

\begin{itemize}

\bigbreak
\item[$\circ$] If $A= \toba(V)$, when is $E(A)$ also a Nichols algebra?

\bigbreak
\item[$\circ$] Is the graded algebra associated to the filtration given by the PBW-basis
on the Lyndon words of $A= \toba(V)$, also a Nichols algebra?

\end{itemize}

\bigbreak
Incidentally, it seems that the determination of the space of point modules
for $U_q^+(\mathfrak g)$ has not been explicited in the literature.

\subsubsection{Invariants}
A natural question in noncommutative geometry is the study of spaces
of invariants under group (or Hopf algebra) actions.
We believe that a suitable setting to discuss it is when the
noncommutative space corresponds to a
Nichols algebra $\toba(V)$  with suitable properties.
Specifically, we propose to study the subalgebra of invariants of
$\toba(V)$ under a coaction of a Hopf algebra $H$.

\bigbreak The first step is to find a good amount of Hopf algebras
$H$ such that $\toba(V)$ is an $H$-comodule algebra. It is
well-known that $\toba(V)$ is a comodule braided Hopf algebra over
the FRT-Hopf algebra $H(c)$ associated to the braided vector space
$(V, c)$, see for example \cite{T}; thus, $\toba(V)$ is a comodule
algebra over any Hopf algebra quotient of $H(c)$. See \cite{Mu2}
for the classification of finite-dimensional quotients of
$\Cc_q[G]$, $G$ a simple algebraic group.

\bigbreak
We stress that these comodule algebra structures are "linear",
that is, they preserve also the comultiplication of $\toba(V)$;
many other structures may arise.
Assume for example that $(V, c) = (V, \tau)$, where $\tau$
is the usual transposition. Then the
Nichols algebra $\toba(V)$ is the symmetric algebra $S(V)$.
The automorphism group of a polynomial algebra is much larger than
the group of linear automorphisms, and the determination
of the former is a classical open problem.

\bigbreak
In the quantum case the situation is much more rigid.
Indeed, assume that $q$ is not a root of 1.
Then $\Aut_{\rm alg} U_q^+(\mathfrak g)$ coincides with
$\Aut_{\rm \ Hopf\ alg} U_q^+(\mathfrak g) \simeq (T \rtimes
\Aut \Delta)$, where $T$ is a maximal torus
and $\Delta$ is the Dynkin diagram, if $\g$ is of type $A_2$ \cite{AlD} or of type $B_2$
\cite{AD}, and conjecturally for all the types.
Again, one is tempted to ask for the class of braided vector spaces $(V, c)$ such that
 $\Aut_{\rm alg} \toba(V) = \Aut_{\rm \ Hopf\ alg} \toba(V)$.

\subsection*{Appendix}

For illustration, we collect some information about finite dimensional
Nichols algebras of rack type. Below we consider braided vector spaces
of rack type, with $q_{ij} \equiv -1$, see \eqref{rack}. The rank of a
Nichols algebra $\toba(V)$ is the dimension of $V$.

\bigbreak
An affine rack is a rack $(A, g)$ where $A$ is a finite abelian group
and $g\in \Aut A$; then $a\trid b = g(b) + (\id - g) (a)$. The
first  four  racks listed below are affine. A subset of a group stable under conjugation
is a rack; so is the set of transpositions in $\mathbb S_n$.
The action for the rack of faces of the cube can be described either
geometrically or as an extension.

\bigbreak
There is no problem to find the space of relations in degree 2; it is
the kernel of $c +\id$. Relations in higher degree (not coming from those in degree 2)
are more difficult to find, as said. For affine racks, a first step is given
in \cite[6.13]{AG2}. Typical relations in degree 4 and 6 are respectively
$$x_0x_1x_0x_1 + x_1x_0x_1x_0 = 0,
\qquad x_0x_1x_2x_0x_1x_2 + x_2x_0x_1x_2x_0x_1 + x_1x_2x_0x_1x_2x_0  = 0.$$
These relations depend upon the order of $-g$.
Most of the computations were done with help of a computer program.
See \cite{G-zoo} for details. No explanation of
the numbers appearing in the table is available until now, but
there are some evident patterns.

\bigbreak
Except for the racks of transpositions in $\mathbb S_4$ and
faces of the cube, all the other racks are simple (they do not project
properly onto a non-trivial rack).
Those two racks are extensions with the same base and fiber but they are not isomorphic.
The similarities between the corresponding Nichols algebras are
explained by a kind of Fourier transform, see \cite[Ch. 5]{AG2}.

\bigbreak More examples of finite-dimensional Nichols algebras of
rack type are given in \cite[Prop. 6.8]{AG2}; they are not of diagonal
type but they arise from Nichols algebras of diagonal type 
by the same kind of Fourier transform.

\begin{table}[t]
\begin{center}
\begin{tabular}{|p{3,7cm}|p{0,5cm}|p{4cm}|p{2,5cm}|p{1,5cm}|}
\hline {\bf Rack} &    {\bf rk}  & {\bf Relations} & $\dim \toba(V)$ &
{\bf top }

\\ \hline  $(\Z/3, \trid^g), g = 2$ \newline (Transpositions in $\mathbb
S_3$)
  & 3  & 5 relations in degree 2     & $12 = 3.2^2$  & $4 = 2^2$

\\ \hline  $(\Z/5, \trid^g), g = 2$
  & 5  & 10 relations in degree 2 \newline  1 relation in degree 4
  & $1280 = 5.4^4$  & $16 = 4^2$

\\ \hline  $(\Z/7, \trid^g), g = 3$
  & 7  & 21 relations in degree 2 \newline  1 relation in degree 6
  & $326592 = 7.6^6$  & $36 = 6^2$

\\ \hline  $(\Z/2 \times \Z/2, \trid^g)$, \newline
$g = \begin{pmatrix} 0 & 1 \\ 1 & 1 \end{pmatrix} $
  & 4  & 8 relations in degree 2 \newline  1 relation in degree 6
  & $72$  & $9 = 3^2$

  \\ \hline  Transpositions in $\mathbb S_4$
  & 6  & 16 relations in degree 2
  & $576$  & $12$

    \\ \hline  Faces of the cube
  & 6  & 16 relations in degree 2
  & $576$  & $12$

\\ \hline  Transpositions in $\mathbb S_5$
  & 10  & 45 relations in degree 2
  & $8294400$  & $40$

\\ \hline
\end{tabular}
\end{center}\end{table}

\subsection*{Acknowledgements}
The author  is grateful to J. Alev, F. Dumas and S. Natale for many
conversations about various aspects of Nichols algebras; several
of the questions in the text arise from discussions with them;
and also to J. T. Stafford and J. Zhang for answers to some
consultations.
The author also  thanks F. Dumas for his warm hospitality during
a visit to the University of Clermont-Ferrand in March 2002
(when this work was began); to S. Catoiu for the kind invitation
to the International Conference in Chicago; and to the IHES, where
this paper was written.


\begin{thebibliography}{AAA}


\bibitem[AA]{AA}
A. Abella and N. Andruskiewitsch, \emph{Compact quantum groups and
the FRT-construction}, Bol. Acad. Ciencias (C\'ordoba) {\bf 63}
(1999), 15-44.

\bibitem[AlD]{AlD}
J. Alev and F. Dumas, \emph{Rigidit\'e des plongements des quotients primitifs minimaux de $U\sb q({\rm
sl}(2))$ dans l'alg\`ebre quantique de Weyl-Hayashi},
Nagoya Math. J. {\bf 143} (1996), 119--146.


\bibitem[AD]{AD}
N. Andruskiewitsch and F. Dumas, \emph{Sur les automorphismes de
$U^+_q(\g)$}, Beitr\"age Algebra Geom., to appear, preprint
(2002).

\bibitem[AG1]{AG}
N. Andruskiewitsch and M. Gra\~na, \emph{Braided {H}opf algebras
over non-abelian groups}, Bol. Acad. Ciencias (C\'ordoba) {\bf 63}
(1999), 45-78. Also in {\tt http://arxiv.org/9802074}.

\bibitem[AG2]{AG2}
\bysame, \emph{From racks to pointed Hopf algebras}, Adv. Math.,
to appear. Also in {\tt http://arxiv.org/0202084}.


\bibitem[AS1]{AS1}
N. Andruskiewitsch and H.-J. Schneider, \emph{Lifting of Quantum
Linear Spaces and Pointed Hopf Algebras of order $ p^3$}, J.
Algebra \textbf{209} (1998), 658--691.


\bibitem[AS2]{AS-adv}
\bysame, \emph{Finite quantum groups and Cartan matrices}, Adv.
Math. \textbf{154} (2000), 1--45.


\bibitem[AS3]{AS-cambr}
\bysame, \emph{Pointed Hopf Algebras}, in ``New directions in
Hopf algebras", 1--68, Math. Sci. Res. Inst. Publ. \textbf{43},
Cambridge Univ. Press, Cambridge, 2002.

\bibitem[AS4]{AS4}
\bysame, \emph{Finite quantum groups over abelian groups of prime
exponent},
 Ann. Sci. Ec. Norm. Super. \textbf{35} (2002), 1--26.


\bibitem[AS5]{AS6}
\bysame, {\it A characterization of quantum groups}, {\tt
math.QA/0201095}, 21 pages.

\bibitem[BGS]{BGS} A. Beilinson, V. Ginzburg and W. S\"orgel, {\it
Koszul duality patterns in representation theory}, J. Amer. Math.
Soc. \textbf{9}, (1996),  473--527.






\bibitem[DCK]{deck} C. De Concini and V. G. Kac,
\emph{Representations of quantum groups at roots of 1},
in ``Operator Algebras, Unitary Representations,
Enveloping Algebras, and Invariant Theory'', ed. A.
Connes {\it et al} (2000); Birkh\"auser, 471--506.


\bibitem[FG]{FG} D. Flores de Chela and J. Green, {\it
Quantum symmetric algebras},
Algebr. Represent. Theory\textbf{ 4} (2001), 55-76.


\bibitem[GK]{GK}
Ginzburg, V.; Kumar, S., {\it Cohomology of quantum groups at roots of unity},
Duke Math. J.  \textbf{69} (1993), 179--198.

\bibitem[GL]{GL} K.Goodearl  and T. Lenagan, \emph{Catenarity in quantum algebras},
J. Pure Appl. Algebra \textbf{111} (1996),  123--142.


\bibitem[G\~n1]{G-cont}
M. Gra\~na, \emph{On Nichols algebras of low dimension},
in New trends in Hopf algebra theory (La Falda, 1999),
Contemp. Math.  \textbf{267} (2000), 111--134.

\bibitem[G\~n2]{G-zoo}
\bysame, \emph{Zoo of Nichols algebras of nonabelian group type},
{\tt http://www-math.mit.edu/$\sim$matiasg/zoo.html}.

\bibitem[Gr]{Gr}  J. Green, {\it
Quantum groups, Hall algebras and quantized shuffles}, in Finite
reductive groups (Luminy, 1994), Progr. Math. \textbf{141},
Birkh\"auser, (1997),  273--290.


\bibitem[Gu]{Gu}  D. Gurevich, {\it
Algebraic aspects of the quantum Yang-Baxter equation}, Leningrad
J. Math. \textbf{2},  (1991),  801--828.


\bibitem[Kh1]{Kh1}  V. Kharchenko,
\emph{An Existence Condition for Multilinear Quantum Operations},
J. Algebra {\bf 217} (1999),  188-228.

\bibitem[Kh2]{Kh2} \bysame, {\it Skew primitive elements in Hopf algebras and related identities},
J. Algebra {\bf 238} (2001), 534-559.


\bibitem[Kh3]{Kh3} \bysame, {\it A combinatorial approach to
the quantification of Lie algebras}, Pacific J. Math. {\bf 203} (2002),
191--233.

\bibitem[L\"o]{L} C. L\"ofwall,
\emph{On the subalgebra generated by the one-dimensional elements
in the Yoneda Ext-algebra} in Algebra, algebraic topology and
their interactions (Stockholm, 1983), 291--338, Lecture Notes in
Math. \textbf{1183}, Springer, Berlin, 1986.


\bibitem[L]{L3} G. Lusztig, \emph{Introduction to quantum groups},
Birkh\"auser, 1993.


\bibitem[MS]{MS}
A. Milinski and H-J. Schneider, {\it  Pointed Indecomposable Hopf
Algebras over Coxeter Groups},
in New trends in Hopf algebra theory (La Falda, 1999),
Contemp. Math. {\bf 267} (2000),  215--236.

\bibitem[M\"u1]{Mu} E. M\"uller,  \emph{Some topics on Frobenius-Lusztig
kernels, I},  J. Algebra  \textbf{206} (1998),  624--658.

\bibitem[M\"u2]{Mu2} \bysame,
\emph{Finite subgroups of the quantum general linear group},
Proc. London Math. Soc. (3)  \textbf{81} (2000), 190--210.

\bibitem[N]{N}  W.D. Nichols, \emph{Bialgebras of type one}, Commun.  Algebra
\textbf{6} (1978), 1521--1552.

\bibitem[R]{R}
J.-E. Roos, \emph{Some non-Koszul algebras}, Progr. Math. 172,
Birkhauser, (1999), 385--389.


\bibitem[Ro1]{Ro1} M. Rosso, \emph{Groupes quantiques et algebres de battage quantiques},
C.R.A.S. (Paris) \textbf{320}   (1995),  145--148.

\bibitem[Ro2]{Ro2} \bysame, \emph{Quantum groups and quantum shuffles},
Inventiones Math. \textbf{133}   (1998),  399--416.

\bibitem[Ro3]{Ro3} \bysame, \emph{Lyndon words and Universal R-matrices},
talk at MSRI, October 26, 1999, available at {\tt http://www.msri.org}.



\bibitem[Sm]{s} Smith, S. Paul,
\emph{Some finite-dimensional algebras related to elliptic
curves}, in
 Representation theory of algebras and related topics (Mexico City, 1994),
 CMS Conf. Proc. \textbf{19} (1996),  315--348, Amer. Math. Soc., Providence, RI.

\bibitem[St]{S} J. T. Stafford, \emph{Noncommutative projective geometry},
in Proceedings of the International Congress of Mathematicians,
Beijing 2002, vol. II (2002), 93--103.

\bibitem[SV]{SV} J. T. Stafford and M. Van den Bergh, \emph{Noncommutative curves and noncommutative
surfaces}, Bull. Amer. Math. Soc. \textbf{38} (2001),  171--216.


\bibitem[T]{T}  M. Takeuchi,
\emph{Survey of braided Hopf algebras},
in New trends in Hopf algebra theory (La Falda, 1999),
 Contemp. Math. \textbf{267}   (2000),  301--324.




\bibitem[VdB]{vdb} Van den Bergh, M.,
\emph{Existence theorems for dualizing complexes over non-commutative graded
and filtered rings},
J. Algebra \textbf{195} (1997), 662--679.



\bibitem[Wa]{Wa} M. Wambst,
{\it Complex  de Koszul quantiques},  Ann. Inst. Fourier
(Grenoble) \textbf{ 43} (1993),  1089--1156.

\bibitem[Wo]{W} S. L. Woronowicz,
{\it Differential calculus on compact matrix pseudogroups (quantum
  groups)},  Comm. Math. Phys. {\bf 122} (1989),  125--170.

\end{thebibliography}
\end{document}